\documentclass[12pt]{article}
\usepackage{amsmath,amsthm,amssymb}

\begin{document}

\title{On solving quadratic congruences}

\author{V.~N.~Dumachev\\
Department of Mathematics,  Voronezh Institute\\ of the Ministry of
Internal Affairs of Russia\\
 394065 Voronezh, Russia \\
e-mail: dumv@comch.ru}

\maketitle

\begin{abstract}
The paper proposes a polynomial formula for solution quadratic
congruences in $\mathbb{Z}_p$. This formula gives the correct answer
for quadratic residue and zeroes for quadratic nonresidue. The
general form of the formula for $p=3 \;\rm{mod}\,4$, $p=5
\;\rm{mod}\,8$ and for $p=9 \;\rm{mod}\,16$ are suggested.


\end{abstract}

Consider the problem of finding solutions to quadratic congruences
\begin{equation}
x^2=y \; \rm{mod}\, p \tag {1}
\end{equation}
in a ring of odd prime characteristic $\mathbb{Z}_p$. As was shown
in \cite{1}, the main tools for the analysis of quadratic
congruences are the Legendre symbol and the Jacobi symbol. However,
the use of these symbols allows us to determine whether the free
term of the congruences is a quadratic residue or a quadratic
nonresidue only.

For further work, we need the following obvious assertions.

\textbf{Definition 1.} $ \forall \mathbb{Z}_p: p \pm 1=0\; \text{mod
}2$.

\textbf{Definition 2.} $\forall \mathbb{Z}_p$:  $p=3\;
\text{mod}\,4$ or $p=1\; \rm{mod}\,4$.

\textbf{Definition 3.} $\forall \mathbb{Z}_p$: if $p \neq 3\;
\text{mod }4$ then $p=1\; \rm{mod}\; 8$ or $p=5\; \rm{mod}\; 8$.

\textbf{Definition 4.}  Euler's criterion for quadratic residues:
\[x^{\frac{p-1}{2}}=1\; \rm{mod}\;p.\]

\textbf{Definition 5.}  Euler's criterion for quadratic nonresidues:
\[x^{\frac{p-1}{2}}=-1\; \rm{mod}\;p.\]

Let  $X \in \mathbb{Z}_p$. Denote by $Y_1 \in \mathbb{Z}_p$
quadratic residues, and by $Y_2 \in \mathbb{Z}_p$ quadratic
nonresidues, then  $X \xleftarrow{f} \{0,Y_1,Y_2\}$. We seek
solution quadratic congruences in polynomial form
\[
x=f(y)=\sum \limits_{i=0}^{p-1}a_iy^i,
\]
where $x \in X$. Obviously, due to the finiteness of field, such a
polynomial always exists. Indeed, the following is verified by
direct substitution

\textbf{Ttheorem 1.} {\it $\forall p=3\; \rm{mod}\;4$, the function
\[
x=\pm y^{\frac{p+1}{4}} \;\rm{mod}\;p.
\]
gives correct solve of congruences (1) for $y \in Y_1$ and the wrong
ones for $y \in Y_2$.}

In this case, each answer must be substituted into formula (1) for
verification, or we need to analyze the free term $y$ using the
Legendre symbol.

The goal of this paper is to find a polynomial
\[
f(y)=\left\{\begin{array}{c}
       x, \text{for}\; y \in Y_1, \\
       0, \text{for}\; y \in Y_2.
     \end{array}\right.
\]
In other words, we are looking for a polynomial that yields the
correct  answer for quadratic residues and zeroes for quadratic
nonresidues.

\textbf{Ttheorem 2.} {\it $\forall p=3\; \text{mod }4$, quadratic
congruences (1) has a solving}
\begin{equation}
x=\pm \frac{p+1}{2}y^{\frac{p+1}{4}}\left(1+y^{\frac{p-1}{2}}\right)
\rm{mod}\;p. \tag{2}
\end{equation}

$\square$ By squaring both sides of the last equation, we get
Euler's criterion for quadratic residues
$x^{\frac{p+1}{2}}=x\;\rm{mod}\; p$. $\blacksquare$

Note that this theorem gives a polynomial that is minimal in the
number of operations. Using the binomial formula, it is easy to
estimate that in just $\mathbb{Z}_p$ we can get $2^{\frac{p-1}{2}}$
different polynomials that give solving quadratic congruences with
zero output for nonresidues. As an example, below given the
polynomials of solution quadratic congruences (1) for small $p=3 \;
\rm{mod}\; 4$.

\[
  \begin{array}{l|l || l|l}
    \mathbb{Z}_p   & \phantom{xxx} x                   &  \mathbb{Z}_p   & \phantom{xxx} x \\ \hline
    \mathbb{Z}_{7} & \pm 4y^2(1+y^3)       &  \mathbb{Z}_{47}& \pm 24y^{12}(1+y^{23}) \\
    \mathbb{Z}_{11}& \pm 6y^3(1+y^5)       &  \mathbb{Z}_{59}& \pm 30y^{15}(1+y^{29}) \\
    \mathbb{Z}_{19}& \pm 10y^5(1+y^9)      &  \mathbb{Z}_{67}& \pm 34y^{17}(1+y^{33}) \\
    \mathbb{Z}_{23}& \pm 12y^6(1+y^{11})   &  \mathbb{Z}_{71}& \pm 36y^{18}(1+y^{35}) \\
    \mathbb{Z}_{31}& \pm 16y^8(1+y^{15})   &  \mathbb{Z}_{79}& \pm 40y^{20}(1+y^{39}) \\
    \mathbb{Z}_{43}& \pm 22y^{11}(1+y^{21})&  \mathbb{Z}_{83}& \pm 42y^{21}(1+y^{41}) \\
  \end{array}
\]

\bigskip
Next, we will look for polynomials generating solutions to quadratic
congruences (1) for $p=1\; \rm{mod}\;4$. This case contains two
options: $p=5\; \rm{mod}\;8$ and $p=1\; \rm{mod}\;8$.

\textbf{Ttheorem 3.} {\it $\forall p=5\; \rm{mod}\;8$, quadratic
congruences (1) has a solution}
\begin{equation}
x=\pm b
y^{\frac{p+3}{8}}\left(1+y^{\frac{p-1}{2}}\right)\left(a+y^{\frac{p-1}{4}}\right)\text{mod
}p, \tag{3}
\end{equation}
where $a^2 =-1 \; \text{mod }p$, $8ab^2 =1 \; \text{mod }p$.

$\square$ Obviously, the factor  $\left(1+y^{\frac{p-1}{2}}\right)$
is necessary in order to vanished all quadratic nonresidues. Will
seek the solution in the form
\[f(y)=g(y)\left(1+y^{\frac{p-1}{2}}\right)\left(a+y^{\frac{p-1}{4}}\right).
\]
By squaring both sides this expression, we get
\[
y=g^2(y)4ay^{\frac{p-1}{4}}\left(1+y^{\frac{p-1}{2}}\right).
\]
If we put $g^2(y)=2y^{\frac{p+3}{4}}$, then the last expression
reduces to the Euler criterion for a quadratic residue.
$\blacksquare$

As an example, below given the polynomials of the solution to
quadratic congruences (1) for small $p=5\; \text{mod }8$.

\[
  \begin{array}{l|c || l|c}
    \mathbb{Z}_p   &   x                   &  \mathbb{Z}_p   &   x \\ \hline
&&&  \\
 \mathbb{Z}_{5} &\left\{ \begin{array}{l}
                  \pm y(1+y^2)(2+y)\\
                  \pm 3y(1+y^2)(3+y)
                  \end{array}\right.
 &  \mathbb{Z}_{37}& \left\{ \begin{array}{l}
                  \pm 8 y^5(1+y^{18})(6+y^9)\\
                  \pm 11y^5(1+y^{18})(31+y^9)
                  \end{array}\right. \\
&&&                   \\
 \mathbb{Z}_{13} &\left\{ \begin{array}{l}
                  \pm y^2(1+y^6)(5+y^3)\\
                  \pm 8y^2(1+y^6)(8+y^3)
                  \end{array}\right.
 &  \mathbb{Z}_{53}& \left\{ \begin{array}{l}
                  \pm 21 y^7(1+y^{26})(23+y^{13})\\
                  \pm 6  y^7(1+y^{26})(30+y^{13})
                  \end{array}\right. \\
&&&                   \\
  \mathbb{Z}_{29} &\left\{ \begin{array}{l}
                  4 y^4(1+y^{14})(17+y^7)\\
                  10y^4(1+y^{14})(12+y^7)
                  \end{array}\right.
& \mathbb{Z}_{61} & \left\{ \begin{array}{l}
                  28 y^8(1+y^{30})(11+y^{15})\\
                  3  y^8(1+y^{30})(50+y^{15})
                  \end{array}\right. \\
  \end{array}
\]

\bigskip
Next, we will look for polynomials generating solution to quadratic
congruences (1) for $p=1\; \rm{mod}\;8$. This case contains two
options: $p=9\; \rm{mod}\;16$ and $p=1\; \rm{mod}\;16$.

\textbf{Ttheorem 4.} {\it $\forall p=9\; \rm{mod}\;16$, quadratic
congruences (1) has a solution}
\[
x=\pm d
y^{\frac{p+7}{16}}\left(1+y^{\frac{p-1}{2}}\right)\left(a+y^{\frac{p-1}{8}}\right)
\left(b+cy^{\frac{p-1}{8}}+y^{\frac{p-1}{4}}\right)\rm{mod}\;p,\tag{4}
\]
where $g_1=2, \; g_2=0, \; g_3=0, \; g_4=0$ or
\[
x=\pm d
y^{\frac{3p+5}{16}}\left(1+y^{\frac{p-1}{2}}\right)\left(a+y^{\frac{p-1}{8}}\right)
\left(b+cy^{\frac{p-1}{8}}+y^{\frac{p-1}{4}}\right)\rm{mod}\;p,
\]
where $g^2_1=0, \; g^2_2=0, \; g^2_3=2, \; g^2_4=0$. Here
\[
\left\{
\begin{array}{l}
  g^2_1=d^2(4 a^2 c+8 a b+4 a c^2+4 b c),  \\
  g^2_2=d^2(2 a^2 b^2+2 a^2+8 a c+4 b+2 c^2),   \\
  g^2_3=d^2(4 a^2 b c+4 a b^2+4 a+4 c), \\
  g^2_4=d^2(4 a^2 b+2 a^2 c^2+8 a b c+2 b^2+2).
\end{array}
\right.\]

$\square$ Obviously, the factor $\left(1+y^{\frac{p-1}{2}}\right)$
is necessary in order to vanished all quadratic nonresidues. Will
seek the solution in the form
\[f(y)=g(y)\left(1+y^{\frac{p-1}{2}}\right)\left(a+y^{\frac{p-1}{8}}\right)\left(b+cy^{\frac{p-1}{8}}+y^{\frac{p-1}{4}}\right).
\]
By squaring both sides this expression, we get 2 cases.

Case 1:
\[
y=g^2(y)y^{\frac{p-1}{8}}\left(1+y^{\frac{p-1}{2}}\right),
\]
where $g(y)=g_1y^{\frac{3p+5}{16}}$.

Case 2:
\[
y=g^2(y)y^{\frac{3p-3}{8}}\left(1+y^{\frac{p-1}{2}}\right),
\]
where $g(y)=g_3y^{\frac{p+7}{16}}$.

Now the last expression reduces to the Euler criterion for a
quadratic residue. $\blacksquare$

As an example, below given the polynomials of the solution to
quadratic congruences (1) for small $p=9\; \rm{mod}\;16$.

\[
  \begin{array}{l|l }
    \mathbb{Z}_p   &   \phantom{xxxxxxxxxxxxxx} x                 \\ \hline
&  \\
 \mathbb{Z}_{41}&\left\{ \begin{array}{l}
                  \pm 5y^3(1+y^{20})(2+y^{5})(30+36y^{5}+y^{10}) \\
                  \pm 15y^3(1+y^{20})(3+y^{5})(9+4y^{5}+y^{10})
                  \end{array}\right. \\
&                   \\
  \mathbb{Z}_{73}& \left\{ \begin{array}{l}
                  \pm 2y^5(1+y^{36})(22+y^{9})(46+3y^{9}+y^{18}) \\
                  \pm 36y^{14}(1+y^{36})(50+y^{9})(8+y^{9}+y^{18})
                  \end{array}\right. \\
&                   \\
 \mathbb{Z}_{89} &\left\{ \begin{array}{l}
                  \pm 12y^6   (1+y^{44})(7+y^{11})(78+80y^{11}+y^{22})\\
                  \pm 37y^{17}(1+y^{44})(3+y^{11})(60+74y^{11}+y^{22})
                  \end{array}\right.  \\
 &                   \\
  \mathbb{Z}_{137}& \left\{ \begin{array}{l}
                  \pm 53y^{ 9}(1+y^{68})(10+y^{17})(100+81y^{17}+y^{34})\\
                  \pm 52y^{26}(1+y^{68})(6 +y^{17})(50+35y^{17}+y^{34})
                  \end{array}\right. \\
  \end{array}
\]

Unfortunately, for the rest $\mathbb{Z}_p$, get a universal
presentation of the polynomial has not yet succeeded. For example,
solutions to quadratic congruences (1)  for $\mathbb{Z}_{17}$ may
have the form
\[
x= \pm y \left(1+y^8\right) \left(y^7+5 y^6+12
y^4+3y^3+15y^2+10y+14\right).
\]

In conclusion, we give simple rules for working with the formulas
obtained. According to Definition 1, for any odd prime $p$, the
upper and lower neighbors are even. It follows from this, for any
prime $p$ either the top or lower neighbor is divided by 4. If the
upper neighbor is divided into 4, then quadratic congruences (1) has
solutions (2). If $(p+1)$ is not divisible by 4, then $(p-1)$ must
divisible by 4. If the result of dividing $\frac{p-1}{4}= \rm{odd}$,
then solutions of quadratic comparison are given by formulas (3). If
$\frac{p-1}{4}= \rm{even}$, then this is equivalent to the fact that
$p=1\; \rm{mod}\; 8$. The following procedure is repeated. If the
result of dividing $\frac{p-1}{8}= \rm{odd}$, then the quadratic
congruences are solved with formulas (4). Condition $\frac{p-1}{8}=
\rm{even}$ is equivalent to $p = 1\; \rm{mod}\; 16$, etc.

\bigskip
\textbf{Examples.} Solve the congruences.

\bigskip
\textbf{1.} $x^2=7\;\rm{mod}\; 103$.

Since $103=3\,\rm{mod}\, 4$ then we will use (2):
\[
x=\pm \frac{103+1}{2}7^{26}(1+7^{51})=\pm 25 \,\rm{mod}\, 103.
\]
It is quadratic residue.

\bigskip
\textbf{2.} $x^2=5\,\rm{mod}\, 107$.

Since $107=3\,\rm{mod}\, 4$ then we will use (2):
\[
x=\pm \frac{107+1}{2}5^{26}(1+5^{51})=0 \,\rm{mod}\, 107.
\]
It is quadratic nonresidue.

\bigskip
\textbf{3.} $x^2=13\,\rm{mod}\, 61$.

Since $61\neq 3\,\rm{mod}\, 4$ then  $61 = 1\,\rm{mod}\, 4$ and we
can not use the formula (2).  Since $61=5\,\rm{mod}\, 8$ then we
will use (3):
\[
x=\pm 28 \cdot 13^{8}(1+13^{30})(11+13^{15})=\pm 47 \,\rm{mod}\, 61.
\]
It is quadratic residue.

\bigskip
\textbf{4.} $x^2=17\,\rm{mod}\, 61$.

Since $61=5\,\rm{mod}\, 8$ then we will use (3):
\[
x=\pm 28 \cdot 17^{8}(1+17^{30})(11+17^{15})=0 \,\rm{mod}\, 61.
\]
It is quadratic nonresidue.

\bigskip
\textbf{5.} $x^2=19\;\rm{mod}\; 137$.

Since $137 \neq 3\,\rm{mod}\, 4$ then  $137 = 1\,\rm{mod}\, 4$ and
we can not use the formula (2).

Since $137 \neq 5\,\rm{mod}\, 8$ then $137 = 1\,\rm{mod}\, 8$ and we
can not use the formula (3).

Since $137=9\,\rm{mod}\, 16$ then we will use (4):
\[
x=\pm 53 \cdot 19^{9}(1+19^{68})(10+19^{17})(100+81\cdot
19^{17}+19^{34})=\pm 108 \,\rm{mod}\, 137.
\]
It is quadratic residue.

\bigskip
\textbf{6.} $x^2=21\,\rm{mod}\, 137$.

Since $137=9\,\rm{mod}\, 16$ then we will use (4):
\[
x=\pm 53 \cdot 21^{9}(1+21^{68})(10+21^{17})(100+81\cdot
21^{17}+21^{34})=0 \,\rm{mod}\, 137.
\]
It is quadratic nonresidue.

\bigskip

\bigskip


\end{document}